\documentclass[a4paper, 12pt]{article}

\usepackage{indentfirst}
\usepackage{graphicx}
\usepackage{epstopdf}
\usepackage{epsfig}
\usepackage{overpic}
\usepackage{overcite}
\usepackage{array}
\usepackage{color}
\usepackage{multirow}
\usepackage{float}
\usepackage{subfigure}
\usepackage{amsmath,amssymb, amsthm}
\usepackage[mathscr]{euscript}
\usepackage{latexsym,bm}
\usepackage{cite}
\usepackage{booktabs}
\usepackage{appendix}

\usepackage{bbm}

\newcommand\comment[1]{}

{\theoremstyle{remark} }

\def\mi{\mathbbm{i}}
\def\me{\mathbbm{e}}
\def\D{\text{d}}

\graphicspath{ {./images/}  }

\begin{document}
\bibliographystyle{unsrt}

\title{Comparison of numerical treatments for the pseudo-differential term in the Wigner equation}
\author{Zhenzhu Chen\footnotemark[2] $^,$\footnotemark[1]}
\renewcommand{\thefootnote}{\fnsymbol{footnote}}
\footnotetext[2]{Institute of Applied Physics and Computational Mathematics, Beijing 100094, China .}
\footnotetext[1]{To
whom correspondence should be addressed. Email:
\texttt{chenzz\_iapcm@163.com}}
\date{\today}
\maketitle

\begin{abstract}
 
Effectively handling the nonlocal pseudo-differential term plays an important role in solving the Wigner equation with high accuracy.
This paper systematically analyzes and compares numerical treatments of the pseudo-differential term under different types of potentials.

\vspace*{4mm}

\end{abstract}

\section{Introduction}
\label{sec:intro}


The Wigner equation has provided a convenient way to render quantum mechanics in phase space since its invention in 1932 \cite{Wigner1932}, and has been successfully applied to semiconductor devices \cite{th:Zhao2000,bk:NedjalkovDimovSelberherr2021} and other fields \cite{TilmaEverittSamsonMunroNemoto2016,WeinbubFerry2018}.   
The pseudo-differential term is the main difference between the Boltzmann equation and the Wigner equation and contains all quantum information \cite{Kontsevich2003}.     
The original form of the pseudo-differential term (denoted as $\Theta_V[f]$) is expressed as follows
\begin{equation}\label{eq:y-form}
 \Theta_V[f](x,k,t) = \frac{1}{2\pi\mi\hbar}\int_{\mathbb{R}}\D y\me^{-\mi ky}D_V(x,y,t)\int_{\mathbb{R}}\D k'\me^{\mi k'y}f(x,k',t)
\end{equation}  
with
\[
   D_V(x,y,t) = V(x+y/2,t) - V(x-y/2,t),
\]
where $f(x,k,t)$ is the Wigner function in phase space $(x,k)$ for the position $x$ and the wave number $k$, $V(x,t)$ is the potential. 
The pseudo-differential term  can be expressed as a convolution form 
\begin{align}
\label{eq:k-form}
 \Theta_V[f](x,k,t) &= \int_{\mathbb{R}}\D k'f(x,k',t)V_w(x,k-k',t),\\
  V_w(x,k,t)&=\frac{1}{2\pi\mi\hbar}\int_{\mathbb{R}}\D y\me^{-\mi ky}D_V(x,y,t),
\end{align}
which is the most commonly used form of the pseudo-differential term. 
The pseudo-differential term can also be characterized  by the Moyal expansion
\begin{equation}\label{eq:m-form}
 \Theta_V[f](x,p,t) = \sum_{l=0}^{+\infty}\frac{(-1)^l\hbar^{2l}}{2^{2l}}\cdot\frac{\nabla_x^{2l+1}V(x,t)}{(2l+1)!}\cdot\nabla_p^{2l+1}f(x,p,t).
\end{equation}
For polynomial potentials, the Wigner equation under the Moyal expansion is  naturally chosen \cite{FurtmaierSucciMendoza2016,ChenXiongShao2019}.
And also the Wigner equation is reformulated using a spectral decomposition of the classical force field instead of the potential \cite{VandeputSoreeMagnus2017}, as shown in Eq.~\eqref{eq:f-form},    
\begin{equation}\label{eq:f-form}
 \Theta_V[f](x,p,t) = \frac{1}{2\pi\hbar}\int_{\mathbb{R}}\D k'\frac{\tilde{F}(k',t)}{k'}\me^{\mi k'x}\left[f(x,k-k'/2,t)-f(x,k+k'/2,t)\right],
\end{equation}
where $F(x,t)=-\nabla_xV(x,t)$ and $\tilde{F}(k,t)=\int_{\mathbb{R}}\D k\me^{-\mi k x}F(x,t)$.

Based on different forms of the pseudo-differential term in the Wigner equation, several different numerical treatments have been proposed. 
This work is the first time to systematically compare and analyze the accuracy and efficiency of numerical treatments for the pseudo-differential term in local potentials, polynomial potentials and general potentials, respectively.
We conclude that the Moyal expansion is the natural choice for polynomial potentials, while numerical treatments for the $y$-integral form \eqref{eq:y-form}, the convolution form \eqref{eq:k-form} and the the force field spectral expansion form \eqref{eq:f-form} all have high accuracy for local potentials.
Meanwhile, we give the prior error estimations $g_{N_\xi}$ of the the $\mathcal{K}$-truncation based on the convolution form under all potentials, so the $\mathcal{K}$-truncation of the pseudo-differential term is our recommended when solving the Wigner equation.   
The typical Gauss barrier, quartic polynomial double-well and discrete potential of resonant tunneling diode (RTD) are further used to verify our conclusions.

The rest of the paper is organized as follows. In Section \ref{sec:2}, numerical treatments for different forms of the pseudo-differential term are given in detail. Section \ref{sec:3} conducts several typical potentials to verify the accuracy of these numerical treatments. The paper ends in Section \ref{sec:conclude} with some conclusions.

\section{Numerical treatments}
\label{sec:2}

This section gives the corresponding numerical treatments according to the different forms of the pseudo-differential term, respectively. 

\begin{itemize}
\item {\bf $\mathcal{Y}$-truncation.} Considering the decay of the density function 
$$\rho(x,y,t):=\frac{1}{2\pi}\int_{\mathbb{R}}\D k'\me^{\mi k'y}f(x,k',t)$$
 when $|y|\rightarrow +\infty$, a direct truncation is applied to the  integral domain in $y$-space, denoted by $\mathcal{Y}=[-L_y/2,L_y/2]$. Further using the decay of the Wigner function in the $k$-space and the  Poisson summation formula, we obtain 
\[
  \frac{1}{2\pi}\int_{\mathbb{R}}\D k'\me^{\mi k'y}f(x,k',t)\approx \frac{\Delta k}{2\pi}\sum_{\mu=-N_\mu}^{N_\mu}f(x,k_\mu,t)\me^{\mi k_{\mu}y},
\]
where $k_{\mu}=\mu\Delta k$ with $\Delta k$ being the spacing in $k$-space. Moreover, in order to maintain the mass conservation, $L_y$ and $\Delta k$ satisfy the condition $L_y\cdot \Delta k = 2\pi$.
 Then, we have the $\mathcal{Y}$-truncation of the pseudo-differential term Eq.~\eqref{eq:y-form} as follows,
\begin{equation}\label{eq:y-truncation}
g^{YT}(x,k,t):=\frac{\Delta k}{2\pi\mi\hbar}\sum_{\mu=-N_\mu}^{N_\mu}f(x,k_\mu,t)\int_{\mathcal{Y}}\D y\me^{-\mi(k-k_{\mu})y}D_V(x,y,t).
\end{equation}
In this work, the Gauss quadrature formula is further used to calculate the integral term in Eq.~\eqref{eq:y-truncation}.
There are two pending parameters in the $\mathcal{Y}$-truncation, namely $L_y$ and $N_\mu$.

\item {\bf $\mathcal{K}$-truncation.}
Since the decay of the Wigner function in the $k$-space, the convolution form Eq.~\eqref{eq:k-form} of the pseudo-differential term can be truncated to integrate over a sufficiently large $k$-domain, denoted by $\mathcal{K}=[-L_k/2,L_k/2]$. Further, using the Poisson summation formula, the convolution kernel $V_w(x,k,t)$ is approximately 
\[
   V_w(x,k,t)\approx \frac{\Delta y}{2\pi\mi\hbar}\sum_{\xi=-N_\xi}^{N_\xi}D_V(x,y_\xi,t)\me^{-\mi k y_\xi},
\]
where $y_\xi =\xi \Delta y$ with $\Delta y$ being the spacing in $y$-space. 
Similarly, in order to ensure mass conservation \cite{XiongChenShao2016}, it is necessary to satisfy $ L_k\cdot \Delta y = 2\pi$.
Therefore, we obtain the $\mathcal{K}$-truncation of the pseudo-differential term according to Eq.~\eqref{eq:k-form}:
\begin{equation}\label{eq:k-truncation}
g^{KT}(x,k,t):=\frac{\Delta y}{2\pi\mi\hbar}\sum_{\xi=-N_\xi}^{+N_\xi}D_V(x,y_\xi,t)\int_{\mathcal{K}}\D k'\me^{-\mi(k-k')y_\xi}f(x,k',t).
\end{equation}
Similarly, the Gauss quadrature formula is applied to calculate the integral term in Eq.~\eqref{eq:k-truncation} and the two pending parameters in the $\mathcal{K}$-truncation are $L_k$ and $N_\xi$.

\item {\bf $\mathcal{M}$-truncation.} The Moyal expansion Eq.~\eqref{eq:m-form} of the pseudo-differential term is an infinite summation. Then the natural way is to truncate the number of summation terms, that is the $\mathcal{M}$-truncation:
 \begin{equation}\label{eq:m-truncation}
g^{MT}(x,k,t):=\sum_{l=0}^{P}\frac{(-1)^l}{\hbar 2^{2l}}\cdot\frac{\nabla_x^{2l+1}V(x,t)}{(2l+1)!}\cdot\nabla_k^{2l+1}f(x,k,t).
\end{equation}
Here, $P$ is the only parameter to be determined.

\item {\bf $\mathcal{F}$-truncation.} 
Using the Poisson summation formula, an approximation of the force field spectrum expansion Eq.~\eqref{eq:f-form} of the pseudo-differential term is obtained 
\begin{equation}\label{eq:psf_ft}
   \Theta_V[f](x,k,t)\approx \frac{\Delta k}{2\pi\hbar}\sum_{\nu=-\infty}^{+\infty}\tilde{F}(k'_\nu,t)\me^{\mi k'_\nu x}\cdot\frac{f(x,k-k'_\nu/2,t)-f(x,k+k'_\nu/2,t)}{k'_\nu},
\end{equation}
where $k'_\nu = \nu\Delta k$ with $\Delta k$ beging the spacing in $k$-space. 
Combining 
\[
  \lim_{k'=0}\tilde{F}(k'_\nu,t)\me^{\mi k'_\nu x}\cdot\frac{f(x,k-k'_\nu/2,t)-f(x,k+k'_\nu/2,t)}{k'_\nu}=-\tilde{F}(0,t)\cdot\nabla_k f(x,k,t),
\]
we further derive the $\mathcal{F}$-truncation
\begin{equation}\label{eq:f-truncation}
\begin{split}
g^{FT}(x,k,t):=&-\frac{\Delta k}{2\pi\hbar}\tilde{F}(0,t)\cdot\nabla_k f(x,k,t)\\
&-\frac{\Delta k}{2\pi\hbar}\sum_{\nu=-N_\nu,\nu\neq 0}^{N_\nu}\tilde{F}(k'_\nu,t)\me^{\mi k'_\nu x}\cdot\frac{f(x,k+k'_\nu/2,t)-f(x,k-k'_\nu/2,t)}{k'_\nu}.
\end{split}
\end{equation}
Here $\tilde{F}(k,t)$ can be approximated by $\int_{\mathcal{X}}\D x F(x,t)\me^{-\mi xk}$ with $\mathcal{X}=[-L_x/2,L_x/2]$ being the domain in $x$-space. $\Delta k$ and $L_x$ satisfy $ L_x\cdot \Delta k = 2\pi$ and the pending parameter in the $\mathcal{F}$-truncation is $N_\nu$.
\end{itemize}

\section{Numerical experiments}
\label{sec:3}

In this section, two typical potentials, local Gauss barrier $V_{loc}(x)=\me^{-x^2/2}$ and unbounded double-well $V_{pol}(x)=(x^2-4)^2$ are employed to test the performance of numerical treatments for the pseudo-differential term. 
 A Gauss wave packet ($f(x,k)=\exp(-x^2/4-4k^2)/\pi$) is used as the Wigner function and the atomic units $\hbar = m = 1$ are applied if not specified.
The numerical performance is evaluated by the $L^\infty$-error $\epsilon_\infty$, defined as follows:
\[
   \epsilon_\infty = \max_{(x,k)\in\mathcal{X}\times\mathcal{K}}\{|\Theta^{\text{ref}}_V[f](x,k) - \Theta^{\text{num}}_V[f](x,k)|\},
\]    
where $\Theta^{\text{ref}}_V[f]$ and $\Theta^{\text{num}}_V[f]$ denote the reference and numerical solution.
 And parameter $g_{N_\xi}$ is defined below to specify the selection criteria for the truncation parameter.
\begin{equation}
g_{N_\xi}=\max_{(x,k)\in \mathcal{X}\times\mathcal{K}}\Delta y\cdot D_V(x,y_{N_\xi},t)\cdot\int_{\mathcal{K}}\D k'\me^{-\mi(k-k')y_{N_\xi}}f(x,k',t).
\end{equation}

\subsection{Local potential -- Gauss barrier}

The analytical formula of the pseudo-differential term corresponding to the Gauss barrier is 
\begin{equation}\label{eq:gb_exact}
  \Theta_V[f](x,k)=\frac{4}{\hbar\pi\sqrt{2\pi}}\int_{\mathbb{R}}\D k'\me^{-2(k-k')^2}\sin(2(k-k')x)\exp(-x^2/4-4k'^2).
\end{equation}
We set $\mathcal{X}=[-10,10]$, $\mathcal{K}=[-2\pi,2\pi]$ and a high-precision Gauss quadrature formula is used to calculate the integral Eq.~\eqref{eq:gb_exact}, which provides the reference solution.
  
 Table~\ref{tab:1} clearly shows that the $\epsilon_\infty$ of the $\mathcal{Y}$-truncation and the value of the potential at $L_y/2$ are of the same order of magnitude. Therefore, for a local potential, the criterion for choosing $L_y$ is to make the value of the potential at $L_y/2$ be zero in the sense of machine precision.
Because of the decay of the Wigner function in the $k$-space, the natural way to choose $N_\mu$ is that $f(x,k)$ is almost zero outside of $[-\Delta k\cdot N_\mu, \Delta k\cdot N_\mu]$. 
Table~\ref{tab:1} also shows that the error reaches the machine precision under $N_\mu=40, L_y=40$. 

\begin{table}[h]
\centering
\caption{\small Gauss barrier. The $L^\infty$-error $\epsilon_\infty$ and $V(L_y/2)$ of $\mathcal{Y}$-truncation \eqref{eq:y-truncation} as a function of $L_y$.}
\label{tab:1}
\begin{tabular}{c||c|c|c|c|c|c}\hline
$L_y$ & 20 & 24 &28 &32 &36 &40\\ 
\hline
$\epsilon_\infty$ & 2.8064e-06 &2.6206e-08&9.3349e-11&1.9501e-13&7.0777e-16 &7.2164e-16\\
\hline
$V(L_y/2)$ & 3.7267e-06 &	1.5230e-08	& 2.2897e-11 & 1.2664e-14	&2.5768e-18		&1.9287e-22\\
\hline
\end{tabular}
\end{table}
Similarly for the $\mathcal{K}$-truncation (Eq.~\eqref{eq:k-truncation}), 
the criterion for choosing $L_k$ is that $f(x,k)$ decays to 0 outside $[-L_k/2,L_k/2]$. 
Table~\ref{tab:2} displays that the $\epsilon_{\infty}$ of the $\mathcal{K}$-truncation and $g_{N_\xi}$ have the same order of magnitude as $N_\xi$. 
Therefore, we choose $N_\xi$ according to the value of $g_{N_\xi}$. And Table~\ref{tab:2} also shows that the error reaches the machine precision $10^{-16}$ under $L_k = 4\pi, N_\xi = 40$. 

\begin{table}[h]
\centering
\caption{\small Gauss barrier. The $L^\infty$-error $\epsilon_\infty$ and $g_{N_\xi}$ of $\mathcal{K}$-truncation \eqref{eq:k-truncation} as a function of $N_\xi$.}
\label{tab:2}
\begin{tabular}{c||c|c|c|c|c|c}\hline
$N_\xi$ & 20 & 24 &28 &32 &36 &40\\ 
\hline
$\epsilon_\infty$ &6.5964e-07& 4.2869e-09&1.7258e-11&2.3092e-14&9.5125e-17&8.9999e-17\\
\hline
 $g_{N_\xi}$ & 6.7178e-07 & 6.8403e-09 & 3.0349e-11 & 5.8859e-14 & 4.7584e-17 &  3.9919e-18 \\
\hline
\end{tabular}
\end{table}
   
Parameter to be determined in the $\mathcal{F}$-truncation (Eq.~\eqref{eq:f-truncation}) is $N_\nu$. 
Table~3 gives the positive correlation between $\epsilon_\infty$ and $g_{N_\nu}:=\Delta k\cdot\tilde{F}_{N_\nu}/k_{N_\nu}$. Therefore, we choose $N_\nu$ so that $g_{N_\nu}$ reaches accuracy to ensure the accuracy and efficiency of $\mathcal{F}$-truncation. Likewise, the error can reach the machine precision at $N_\nu=30$.
\begin{table}[h]
\centering
\caption{\small Gauss barrier. The $L^\infty$-error $\epsilon_\infty$ and $g_{N_\nu}$ of $\mathcal{F}$-truncation \eqref{eq:f-truncation} as a function of $N_\nu$. }
\label{tab:3}
\begin{tabular}{c||c|c|c|c|c|c}\hline
$N_\nu$ & 10 & 14 &18 &22 &26 &30\\ 
\hline
$\epsilon_\infty$ & 2.5228e-04 & 1.4269e-06 & 1.5855e-09 & 1.3313e-11 &
7.3552e-16 & 7.3552e-16\\
\hline
 $g_{N_\nu}$ & 0.0023 & 1.9795e-05 & 3.5754e-08 & 3.2367e-13 & 1.0219e-15 &  1.6171e-20 \\
\hline
\end{tabular}
\end{table}

\begin{figure}[ht!]
\includegraphics[height=0.25\textwidth,width = 0.3\textwidth]{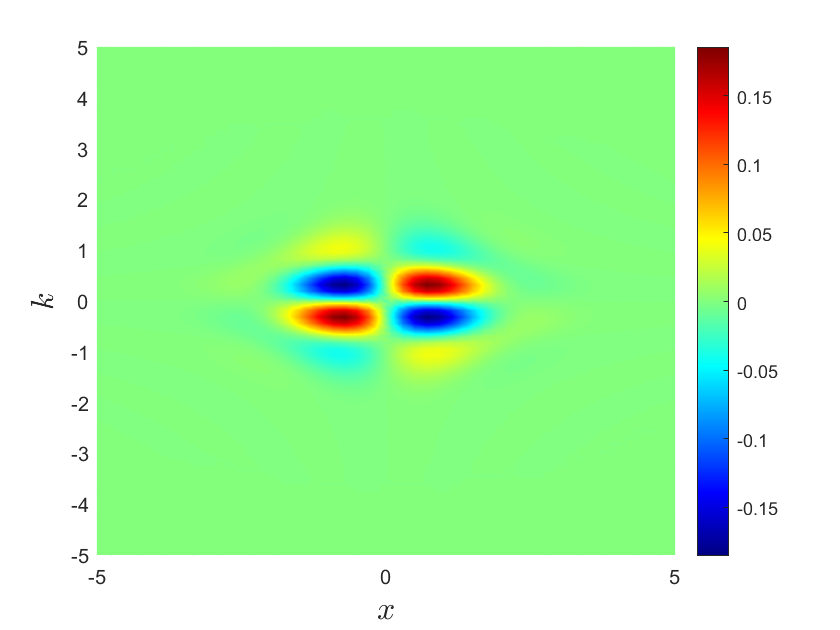}
\includegraphics[height=0.25\textwidth,width = 0.3\textwidth]{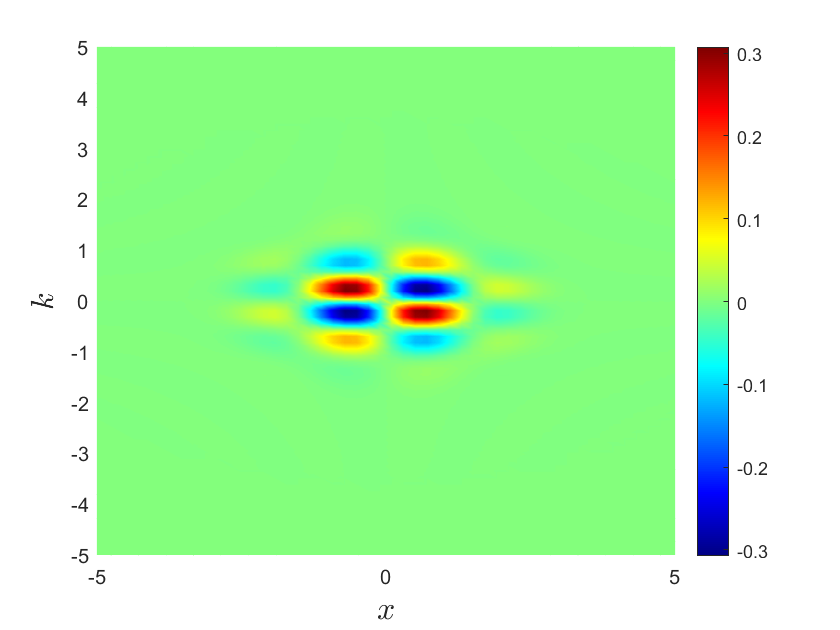}
\includegraphics[height=0.25\textwidth,width = 0.3\textwidth]{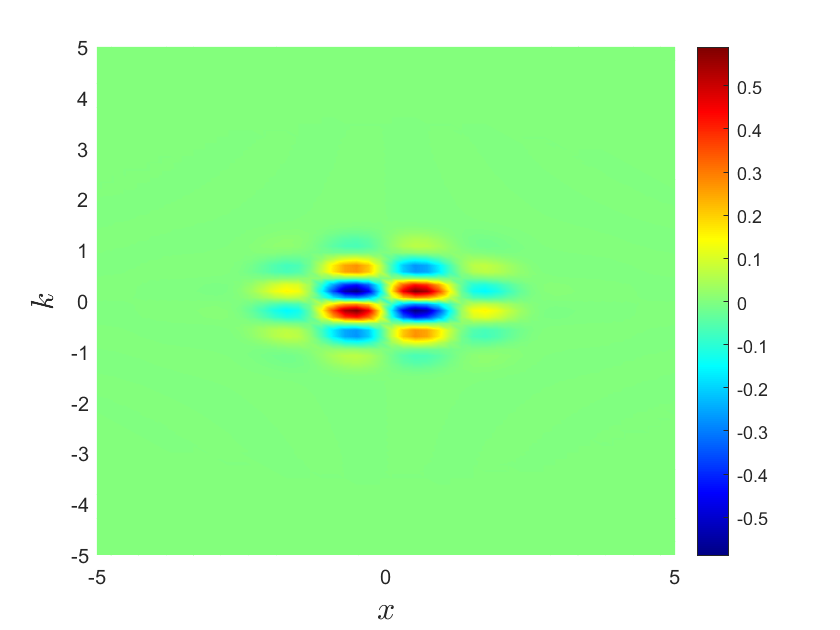}

\includegraphics[height=0.25\textwidth,width = 0.3\textwidth]{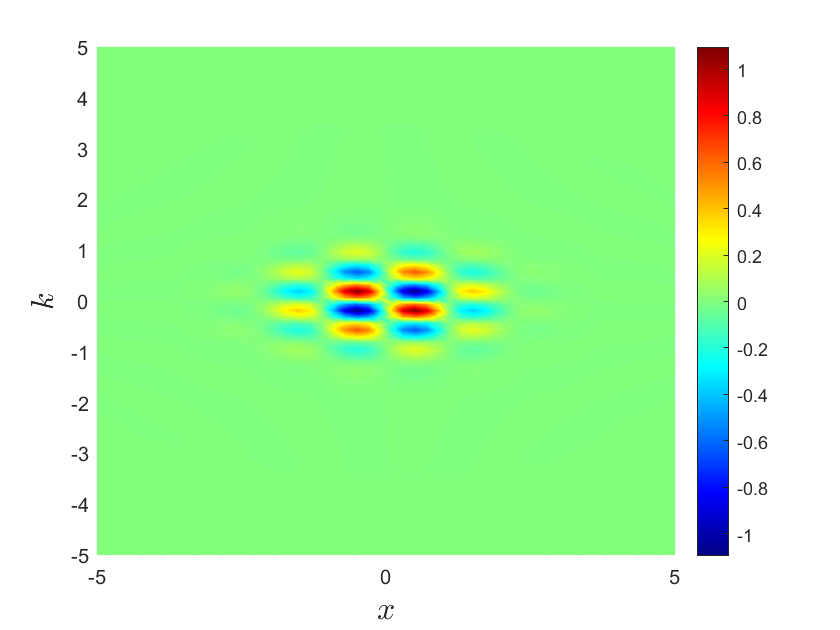}
\includegraphics[height=0.25\textwidth,width = 0.3\textwidth]{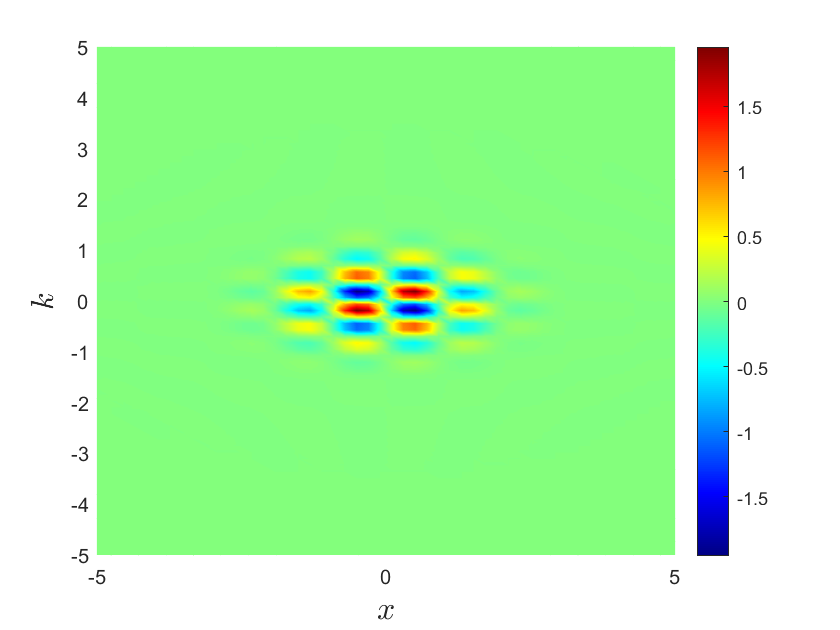}
\includegraphics[height=0.25\textwidth,width = 0.3\textwidth]{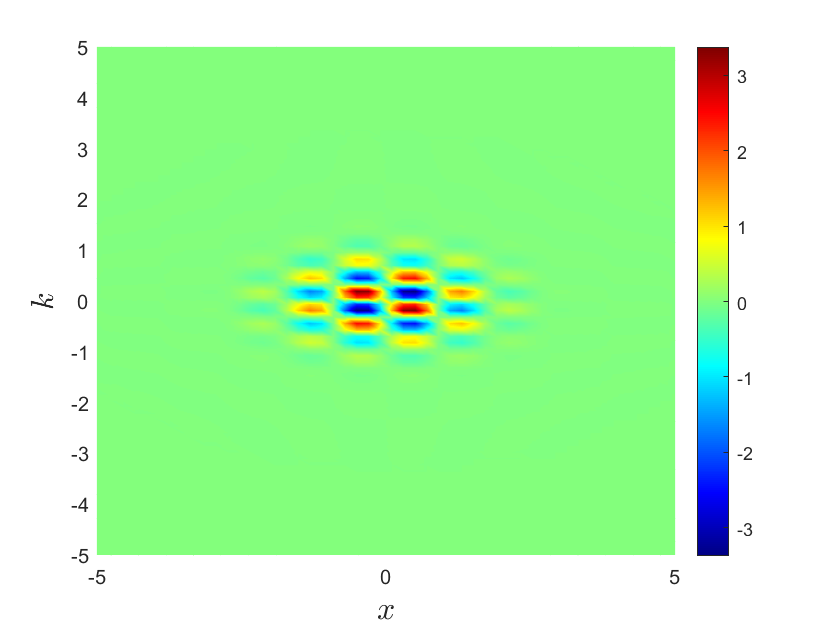}
 \caption{\small Gauss barrier. Distributions of the $\mathcal{M}$-truncation \eqref{eq:m-truncation} with $l=0,1,2,3,4,5$. It is clear that there is no trend of convergence.}
 \label{fig:g_error_m}
\end{figure}
 Fig.~\ref{fig:g_error_m} displays the error distribution under $P=0,1,2,3,4,5$ respectively. It is clearly observed that the error does not decrease with the increase of $P$. This is because the value of high-order derivatives of the Gauss barrier and the Gauss wave packet are both important and cannot be ignored.

\subsection{Unbounded polynomial potential -- Double-well}

The analytical formula of the pseudo-differential term corresponding to the double-well is 
\begin{equation}\label{eq:db_exact}
  \Theta_V[f](x,k)= 4a(x^3-b^2x)\cdot\nabla_kf(x,k)/\hbar-ax\cdot\nabla_k^3f(x,k)/\hbar,
\end{equation}
which is the Moyal expansion. Therefore, the $\mathcal{M}$-truncation is exact for such polynomial potential. We set $a=1$, $b=2$, $\mathcal{K}=[-2\pi,2\pi]$ and $\mathcal{X}=[-15,15]$.

The criterion for choosing $L_y$ for local potentials is given in the previous section, but this does not apply to the unbounded potential. Table~\ref{tab:4} shows that the $L^\infty$-error $\epsilon_\infty$ first decreases with the increase of $L_y$ and finally remains in the order of $10^{-12}$. Here, the appropriate $L_y$ can only be selected through multiple numerical experimental comparisons and the way choosing $N_\mu$  is also to let the Wigner function be zero outside of domain $[-\Delta k\cdot N_\mu,\Delta k\cdot N_\mu]$.  

\begin{table}[h]
\centering
\caption{\small Double-well. The $L^\infty$-error $\epsilon_\infty$ of $\mathcal{Y}$-truncation \eqref{eq:y-truncation} as a function of $L_y$.}
\label{tab:4}
\begin{tabular}{c||c|c|c|c|c|c}\hline
$L_y$ & 30 &35 &40 &45 &50 &55\\ 
\hline
$\epsilon_\infty$ &1.0237e-04& 8.6500e-07 &3.2770e-09& 5.8891e-12& 1.1939e-12&3.8723e-12\\ 
\hline
\end{tabular}
\end{table}
Similarly, the criterion for choosing $L_k$ is that the Wigner function $f(x,k)$ decays to 0 outside $[-L_k/2,L_k/2]$. 
Table~\ref{tab:5} also displays that the $\epsilon_{\infty}$ and $g_{N_\xi}$ have the same order of magnitude as $N_\xi$ for an unbounded polynomial potential. However, it clearly shows that the error does not decrease with the increase of $N_\xi$ after reaching $10^{-11}$. There is an error in the $\mathcal{K}$-truncation because the potential does not decay to 0 at the boundary. Therefore, for the unbounded potential, the criterion for choosing $N_\xi$ by the value of $g_{N_\xi}$ is also feasible. 

\begin{table}[h]
\centering
\caption{\small Double-well. The $L^\infty$-error $\epsilon_\infty$ and $g_{N_\xi}$ of $\mathcal{K}$-truncation \eqref{eq:k-truncation} as a function of $N_\xi$.}
\label{tab:5}
\begin{tabular}{c||c|c|c|c|c|c}\hline
$N_\xi$ & 20 & 30 &40 &50 &60 & 70 \\ 
\hline
$\epsilon_\infty$ &0.0896 & 6.5680e-05 &1.7387e-09& 1.4153e-11& 2.2177e-11&7.1682e-11\\
\hline
 $g_{N_\xi}$ & 0.2155 & 3.0608e-04 & 1.3133e-08 & 1.1937e-11 & 1.9085e-11 & 2.7313e-11 \\
\hline
\end{tabular}
\end{table}

Fixed $\Delta k = 2\pi/L_x$ of the $\mathcal{F}$-truncation, we find that the $L^\infty$-error $\epsilon_\infty$ of the unbounded polynomial potential first decreases slowly with the increase of $N_\nu$, then quickly drops to $10^{-9}$ when $N_\nu=80$, and finally stays at $10^{-10}$ and no longer decreases.
Since the polynomial potential is not attenuated, there is an error in the calculation of $\tilde{F}(k)$ in a bounded domain $\mathcal{X}$. Therefore, the truncated error will not drop to $10^{-16}$. 
\begin{table}[h]
\centering
\caption{\small Double-well. The $L^\infty$-error $\epsilon_\infty$ of $\mathcal{F}$-truncation \eqref{eq:f-truncation} as a function of $N_\nu$.}
\label{tab:6}
\begin{tabular}{c||c|c|c|c|c}\hline
$N_\nu$ & 30 & 40 &50 &60 &70 \\ 
\hline
$\epsilon_\infty$ & 5.9614 & 3.7693 & 2.5220 & 1.6360 & 0.0078 \\
\hline
\hline
$N_\nu$ & 80 & 90 &100 & 110& 120 \\ 
\hline
$\epsilon_\infty$ &6.3389e-09&1.8493e-10&1.8493e-10&1.8493e-10&1.8493e-10\\
\hline
\end{tabular}
\end{table}

According to the numerical analysis and comparison in Sections 3.1 and 3.2, numercial treatments based on $\mathcal{Y}$-truncation, $\mathcal{K}$-truncation and $\mathcal{F}$-truncation are all of high accuracy for local potentials, while for polynomial potentials, the $\mathcal{M}$-truncation is exact. At the same time, for both local potentials and polynomial potentials, we give the prior error estimation of the $\mathcal{K}$-truncation, namely $g_{N_\xi}$.

\subsection{General potentials}

There is a type of potential, which has no specific expression and only   
gives the values at some collocation points, such as the potential determined by solving the Poisson equation in semiconductor device RTD.   
Since the parameter selection criteria and accuracy of the $\mathcal{K}$-truncation are both given under the local potential and unbounded potential in Sections 3.1 and 3.2, we recommend to use the $\mathcal{K}$-truncation for this general potential. Here, we set domain parameters as $\mathcal{K}=[-2\pi,2\pi]$, $\mathcal{X}=[-10,10]$.

\begin{table}[h]
\centering
\caption{\small General potential. The $L^\infty$-error $\epsilon_\infty$ and $g_{N_\xi}$ of $\mathcal{K}$-truncation \eqref{eq:k-truncation} as a function of $N_\xi$.}
\label{tab:7}
\begin{tabular}{c||c|c|c|c|c|c}\hline
$N_\xi$ & 20 & 30 &40 &45 &50 & 60 \\ 
\hline
$\epsilon_\infty$ &3.9148e-05&5.0742e-09& 5.7193e-14&9.0553e-16&
8.6129e-16&6.5312e-16 \\
\hline
 $g_{N_\xi}$ &1.4844e-04& 2.7443e-08&4.8822e-13& 5.9262e-16& 1.1688e-16&2.4655e-16\\
\hline
\end{tabular}
\end{table}
Table~\ref{tab:7} gives the $L^\infty$-error $\epsilon_\infty$ and $g_{N_\xi}$ of the $\mathcal{K}$-truncation, which have almost the same magnitude as we stated before. It further verified our criterion for choosing the parameter of $\mathcal{K}$-truncation. Furthermore, because the  $L^\infty$-error $\epsilon_\infty$ has the same order as $g_{N_\xi}$, $g_{N_\xi}$ also provides a prior error estimation for the $\mathcal{K}$-truncation.

\section{Conclusion}
\label{sec:conclude}

In this paper, we analyze and compare the accuracy of numerical treatments for the pseudo-differential term in local potentials and unbounded polynomial potentials, respectively. 
 For a local potential, $\mathcal{Y}$-truncation, $\mathcal{K}$-truncation and $\mathcal{F}$-truncation are almost machine-accurate. And for an unbounded polynomial potential, $\mathcal{M}$-truncation is exact, while $\mathcal{Y}$-truncation and $\mathcal{K}$-truncation also have high accuracy and efficiency. Meanwhile, we give a prior error estimation  of $\mathcal{K}$-truncation for general potentials. Therefore, $\mathcal{K}$-truncation is 
recommended.

\section*{Acknowledgement}

This work was supported by China Postdoctoral Science Foundation (No.2021M690467). The author is grateful to the useful discussions with Sihong Shao.


\begin{thebibliography}{10}

\bibitem{Wigner1932}
E.~Wigner.
\newblock On the quantum corrections for thermodynamic equilibrium.
\newblock {\em Phys. Rev.}, 40:749--759, 1932.

\bibitem{th:Zhao2000}
P.~Zhao.
\newblock {\em Wigner-{P}oisson Simulation of Quantum Devices}.
\newblock PhD thesis, Stevens Institute of Technology, 2000.

\bibitem{bk:NedjalkovDimovSelberherr2021}
M.~Nedjalkov, I.~Dimov, and S.~Selberherr.
\newblock {\em Stochastic approaches to electron transport in micro- and
  nanostructures}.
\newblock Birkh\"{a}user, Switzerland, 2021.

\bibitem{TilmaEverittSamsonMunroNemoto2016}
T.~Tilma, J.~Everitt, J.~H. Samson, W.~J. Munro, and K.~Nemoto.
\newblock Wigner function for arbitrary quantum systems.
\newblock {\em Phys. Rev. Lett.}, 117:180401, 2016.

\bibitem{WeinbubFerry2018}
J.~Weinbub and D.~K. Ferry.
\newblock Recent advances in {Wigner} approaches.
\newblock {\em April}, 5:041104, 2018.

\bibitem{Kontsevich2003}
M.~Kontsevich.
\newblock Deformation quantization of {Poisson} manifolds.
\newblock {\em Lett. Math. Phys.}, 66:157--216, 2003.

\bibitem{FurtmaierSucciMendoza2016}
O.~Furtmaier, S.~Succi, and M.~Mendoza.
\newblock Semi-spectral method for the {Wigner} equation.
\newblock {\em J. Comput. Phys.}, 305:1015--1036, 2016.

\bibitem{ChenXiongShao2019}
Z.~Chen, Y.~Xiong, and S.~Shao.
\newblock {Numerical methods for the Wigner equation with unbounded potential}.
\newblock {\em J. Sci. Comput.}, 79:345--368, 2019.

\bibitem{VandeputSoreeMagnus2017}
M.~L.~Van de~Put, B.~Sor\'ee, and W.~Magnus.
\newblock {Efficient solution of the {Wigner-Liouville} equation using a
  spectral decomposition of the force field}.
\newblock {\em J. Comput. Phys.}, 350:314--325, 2017.

\bibitem{XiongChenShao2016}
Y.~Xiong, Z.~Chen, and S.~Shao.
\newblock An advective-spectral-mixed method for time-dependent many-body
  {Wigner} simulations.
\newblock {\em SIAM J. Sci. Comput.}, 38:B491--B520, 2016.

\end{thebibliography}
\end{document}